\documentclass[11pt]{amsart}
\usepackage{amssymb, amstext, amscd, amsmath}
\usepackage{mathtools}

\makeatletter
\def\@cite#1#2{{\m@th\upshape\bfseries%
[{#1\if@tempswa{\m@th\upshape\mdseries, #2}\fi}]}}
\makeatother

\mathtoolsset{centercolon} 

\theoremstyle{plain}
\newtheorem{thm}{Theorem}[section]

\theoremstyle{definition}

\newtheorem{conj}[thm]{Conjecture}

\renewcommand{\labelenumi}{(\roman{enumi}) }

  \newcommand{\A}{{\mathcal{A}}}
  \newcommand{\B}{{\mathcal{B}}}
  
  \newcommand{\D}{{\mathcal{D}}}

\renewcommand{\H}{{\mathcal{H}}}

  \newcommand{\K}{{\mathcal{K}}}
\renewcommand{\L}{{\mathcal{L}}}
  \newcommand{\M}{{\mathcal{M}}}
  \newcommand{\N}{{\mathcal{N}}}
\renewcommand{\O}{{\mathcal{O}}}

  \newcommand{\R}{{\mathcal{R}}}
\renewcommand{\S}{{\mathcal{S}}}
  \newcommand{\T}{{\mathcal{T}}}

\newcommand{\bB}{{\mathbb{B}}}
\newcommand{\bC}{{\mathbb{C}}}
\newcommand{\bD}{{\mathbb{D}}}

\newcommand{\bN}{{\mathbb{N}}}

\newcommand{\bR}{{\mathbb{R}}}
\newcommand{\bT}{{\mathbb{T}}}
\newcommand{\bZ}{{\mathbb{Z}}}

\newcommand{\ep}{\varepsilon}
\renewcommand{\phi}{\varphi}
\newcommand{\upchi}{{\raise.35ex\hbox{\ensuremath{\chi}}}}

\newcommand{\fB}{{\mathfrak{B}}}
\newcommand{\fK}{{\mathfrak{K}}}
\newcommand{\fL}{{\mathfrak{L}}}
\newcommand{\fS}{{\mathfrak{S}}}

\newcommand{\rC}{{\mathrm{C}}}

\newcommand{\foral}{\text{ for all }}
\newcommand{\qand}{\quad\text{and}\quad}
\newcommand{\qfor}{\quad\text{for}\ }
\newcommand{\qforal}{\quad\text{for all}\ }

\newcommand{\Ad}{\operatorname{Ad}}
\newcommand{\Alg}{\operatorname{Alg}}
\newcommand{\Aut}{\operatorname{Aut}}
\newcommand{\diag}{\operatorname{diag}}
\newcommand{\dist}{\operatorname{dist}}

\newcommand{\ind}{\operatorname{ind}}

\newcommand{\Lat}{\operatorname{Lat}}
\newcommand{\Mult}{\operatorname{Mult}}

\newcommand{\rank}{\operatorname{rank}}

\newcommand{\Sp}{\operatorname{sp}}
\newcommand{\Tr}{\operatorname{Tr}}

\newcommand{\ca}{\mathrm{C}^*}

\newcommand{\Fd}{\mathbb{F}_d^+}
\newcommand{\Fock}{\ell^2(\Fd)}

\newcommand{\ip}[1]{\langle #1 \rangle}
\newcommand{\ol}{\overline}
\newcommand{\sot}{\textsc{sot}}
\newcommand{\wot}{\textsc{wot}}

\begin{document}

\title[Arveson's Legacy]{The mathematical legacy of William Arveson}

\author[K. R. Davidson]{Kenneth R. Davidson}
\thanks{Author partially supported by an NSERC grant.}
\address{KENNETH R. DAVIDSON, Department of Pure Mathematics,
University of Waterloo, Waterloo, ON N2L 3G1, CANADA}
\email{krdavids@uwaterloo.ca}

\begin{abstract}
This  is a retrospective of some of William B. Arveson's many contributions to 
operator theory and operator algebras.
\end{abstract}

\dedicatory{Dedicated to the memory of Bill Arveson}



\date{}
\maketitle

\section{Introduction}

William Barnes Arveson, born in November 1934, completed his doctorate in 1964 
at UCLA under the supervision of Henry Dye.
After an instructorship at Harvard, Bill began a long career at the University of
California, Berkeley. He died in November 2011, shortly before his
77th birthday, of complications resulting from surgery.

Arveson's work was deep and insightful, and occasionally completely revolutionary.
When he attacked a problem, he always set it in a general framework,
and built all of the infrastructure needed to understand the workings.
This is one of the reasons that his influence has been so pervasive in many areas
of operator theory and operator algebras.

He worked in both operator theory and operator algebras. I think it is fair to say
that he did not really distinguish between the two areas, and proved over and over
again that these areas are inextricably linked.  In \cite{Arv1967b}, he writes in the 
introduction, and in fact demonstrated many times, that:
\begin{quotation}
``In the study of families of operators on Hilbert space, the self-adjoint algebras 
have occupied a preeminent position.
Nevertheless, many problems in operator theory lead obstinately toward questions 
about algebras that are not necessarily self-adjoint.''
\end{quotation}
Bill's work moves across the landscape of operator theory and operator algebras,
and he developed many foundational ideas which wove the subject 
together as an integrated whole.

It isn't possible in a short article to review all of Bill's contributions. 
I will cover many of the highlights. The choices made reflect my personal taste, 
and I extend my apologies to those who feel that their favourite Arveson paper
is omitted.  There will be one glaring omission---$E_0$-semigroups and CP semigroups.
Beginning with \cite{Arv1989}, Arveson wrote 24 papers (more than 650 pages) 
and a book (450 pages) on this topic.
It is too large a topic to cover in this article, and 
Masaki Izumi has written a survey \cite{Izumi} of this area to appear in this same volume.

Bill Arveson was well-known to be an outstanding expositor. His papers were
always meticulously written. In addition, he wrote three graduate level texts.
His \textit{An invitation to C*-algebras} \cite{Arv1976} is distinguished from other 
C*-algebra texts in that it develops in detail the structure of type I C*-algebras and 
their representation theory.
\textit{A short course on spectral theory} \cite{Arv2002a} is a brief introduction to 
Banach algebras and operator theory.
\textit{Noncommutative dynamics and E-semigroups} \cite{Arv2003} is a full development
of the theory of $E_0$-semigroups.

Bill was also a dedicated supervisor.  According to the Mathematics Genealogy Project, 
he had 27 doctoral students. At the memorial ceremony in Berkeley in February, many
of his former students paid their respects and told stories of Bill as a mentor,
inspiration and friend. He will be greatly missed. But his mathematics will live on.

\section{Dilation theory} \label{S:dilation}

I will begin with Arveson's contributions to dilation theory. 
While this is not his first work, it is  perhaps his most influential. 
It has had ramifications in many aspects of operator algebras and operator theory. 

In 1953, Bela Sz.Nagy \cite{SzN} showed that for every operator $T$ on Hilbert space 
of norm at most one, there is a \textit{dilation} to a unitary operator $U$ on a larger 
Hilbert space $\K = \K_- \oplus \H \oplus \K_+$ of the form
\[ U = \begin{bmatrix} *&0&0 \\ *&T&0 \\ *&*&* \end{bmatrix} ,\]
where  $*$ represents an unspecified operator.
In particular, the map from polynomials in $U$ to polynomials in $T$
obtained by compression to $\H$ is easily seen to be a contractive homomorphism. 
This recovers the von Neumann inequality for contractions:
\[ \| p(T)\| \le \|p\|_\infty := \sup_{|z|\le1} |p(z)| \qforal p \in \bC[z]. \]

When $T$ does not have a unitary summand,
this map extends to yield a map from the algebra $H^\infty$ of bounded analytic
functions on the unit disk into the weakly closed algebra that $T$ generates.
This allows the transference of powerful techniques from function theory to operator theory.
See \cite{SF,SF2} for more in this direction.

These methods extended, in part, to pairs of commuting contractions.
In 1964, Ando \cite{Ando} showed that any commuting pair of contractions 
simultaneously dilate to a pair of commuting unitaries. 
In particular, they also satisfy von Neumann's inequality for polynomials in
two variables, where the norm on polynomials is the supremum over the bidisk $\bD^2$.
In the one variable case, the minimal dilation is unique in a natural sense.
However, difficulties arise in the two variable case because this dilation is no longer unique.

For commuting triples, there is a serious road block. 
Briefly skipping ahead in the history,
Parrott \cite{Par} used Arveson's ideas to construct three commuting contractions
which do not dilate to commuting unitaries. 
Furthermore, examples developed later from the work of Varo\-poulos \cite{Var} on tensor 
product norms showed that even von Neumann's inequality can fail.

In two papers in Acta Math.\ \cite{Arv1969, Arv1972}, Arveson revolutionized
dilation theory by reformulating it in the framework of an arbitrary operator algebra. 
The algebra in Sz.Nagy's case is the disk algebra $A(\bD)$,
and for Ando is the bidisk algebra $A(\bD^2)$.

There were a number of key ideas that are central to Arveson's approach.
First, an operator algebra $\A$ was a considered as a subalgebra of a C*-algebra $\fB$.
Second, the correct representations were \textit{completely contractive},
and yielded \textit{completely positive} maps on $\ol{\A+\A^*}$.
Third, he claimed that every operator algebra $\A$ lived inside a unique canonical
smallest C*-algebra called the \textit{C*-envelope}.
These ideas completely changed the way we look at dilations. 

To elaborate, the $k\times k$ matrices over the C*-algebra $\fB$, $\M_k(\fB)$, is itself
a C*-algebra and thus carries a unique norm. This induces a family of norms on
the matrix algebras over $\A$. Today we say that this gives $\A$ an operator space structure.
A map $\phi$ from $\A$ into $\B(\H)$ (or into any C*-algebra) induces maps $\phi_k$
from $\M_k(\A)$ into $\M_k(\B(\H)) \simeq \B(\H^{(k)})$ by acting by $\phi$ on each entry:
\[ \phi_k( \big[ a_{ij} \big] = \big[ \phi(a_{ij}) \big] .\]
We say that $\phi$ is \textit{completely bounded} when
\[ \|\phi\|_{cb} := \sup_{k\ge1} \|\phi_k\| < \infty .\]
It is \textit{completely contractive} when $\|\phi\|_{cb}\le1$.
When $\A$ is an \textit{operator system} (a unital self-adjoint subspace), 
we say $\phi$ is \textit{completely positive} when $\phi_k$ is positive for all $k \ge1$.

Completely positive maps on C*-algebras were introduced by Stinespring \cite{Stine},
where he proved a fundamental structure theorem for completely positive maps. 
The notion was used by St\o rmer \cite{Storm} in a study of positive maps on C*-algebras. 
However, Arveson really recognized the power of this idea, and 
turned completely positive maps into a major tool in operator theory.  

Every completely contractive unital map $\rho$ of a unital operator algebra $\A$ 
extends uniquely to a completely positive map $\tilde\rho$ on the operator system $\ol{\A+\A^*}$. 
Arveson's Extension Theorem established that this map extends to a completely positive
map on the C*-algebra $\fB$.  

\begin{thm}[Arveson's extension theorem 1969]
Let $\S$ be an operator system in a unital C*-algebra $\fB$.
If $\phi$ is a completely positive map from $\S$ into $\B(\H)$,
then there is a completely positive extension $\Phi$ of $\fB$ into $\B(\H)$.
\end{thm}

In modern terminology, this says that $\B(\H)$ is
\textit{injective} in the category of operator systems with completely positive maps.  
This is a generalization of Krein's theorem that positive linear functionals have
a positive extension. The norm is not increased because $\|\Phi\|=\|\Phi(1)\|=\|\phi(1)\|=\|\phi\|$. 

Stinespring's Theorem shows that $\Phi$ is a corner of a $*$-representation $\pi$ of $\fB$. 
When restricted to $\A$, we obtain a completely contractive
representation of the form
\[
 \pi(a) = \begin{bmatrix} * & 0 & 0 \\ * & \rho(a) & 0 \\ * & * & * \end{bmatrix}
 \qforal a \in \A .
\]
This is called a \textit{$*$-dilation} of $\rho$. So we have

\begin{thm}[Arveson's dilation theorem 1969]
Let $\A$ be a unital operator algebra, and
let $\rho:\A \to \B(\H)$ be a representation. 
The following are equivalent:
\begin{enumerate}
\item $\rho$ is completely contractive.
\item $\tilde\rho$ is completely positive.
\item $\rho$ has a $*$-dilation.
\end{enumerate}
\end{thm}

This result was highlighted in Paulsen's book \cite{Paul1} and the expanded 
revision \cite{Paul2}.
In Arveson's paper, this fundamental result is a remark following Theorem 1.2.9.

Arveson was focussed on the notion of a \textit{boundary representation}, that is,
an irreducible $*$-representation $\pi$  of $\fB$ with the property that $\pi|_\A$ has
a unique completely positive extension to $\fB$. These representations are the
non-commutative analogue of points in the spectrum of a function algebra which
have a unique representing measure---points in the Choquet boundary.
When there are sufficiently many boundary representations to represent $\A$
completely isometrically, then there is a quotient of $\fB$ which is the unique 
minimal C*-algebra containing $\A$ completely isometrically. 
This is the \textit{C*-envelope} of $\A$.

Arveson was only able to establish the existence of this C*-envelope in a variety of
special cases. For example, if $\fB$ contains the compact operators and the quotient
by the compacts is not completely isometric on $\A+\A^*$, then $\fB$ itself is the 
C*-envelope.
He analyzed a variety of singly generated examples, always from the boundary 
representation viewpoint.
In particular, Arveson introduced the notion of \textit{matrix range} of an operator $T$:
\[
 W_n(T) =  
 \big\{ \phi(T) : \phi:\ca(T)\to\M_n \text{  unital completely positive} \big\}  
 \qfor n\ge1.
\]
In \cite{Arv1970}, he proved that the matrix range of an irreducible compact operator
formed a complete unitary invariant. 
Finding unitary invariants for operators was an important open problem.
This provides a new invariant even for matrices.
This was strengthened in \cite{Arv1972} to show that for any irreducible operator $T$ 
such that $\ca(T)$ contains the compact operators,
the matrix range is again a complete unitary invariant.

\begin{thm}[Arveson 1972]
If $T_1$ and $T_2$ are irreducible operators such $\ca(T_i)$
contains the compact operators, then $T_1$ and $T_2$ are unitarily equivalent
if and only if 
\[ W_n(T_1) = W_n(T_2) \qforal n \ge 1 .\]
\end{thm}

A decade later, Hamana \cite{Ham} established the existence of a
unique minimal injective operator system containing $\A$, the \textit{injective envelope},
and from this deduced the existence of the C*-envelope in complete generality.
However it shed little light on the existence of boundary representations.
In 2005, Dritschel and McCullough \cite{DritsMcC} produced a very different proof
of the existence of the C*-envelope based on ideas of Agler \cite{Agler}.
It was more direct, and it showed that every representation had a \textit{maximal dilation},
which is a representation such that any further dilation (to a completely contractive 
representation) can be obtained only by adjoining another representation by a direct sum.
This property in another guise was recognized by Muhly and Solel \cite{MS_bound},
but their result relied on Hamana's theorem. This approach did not deal with the question of
irreducible representations either, so it still did not establish the existence of boundary representations.
Arveson revisited this question in the light of this new proof, and established the existence
of sufficiently many boundary representations when $\A$ is separable \cite{Arv2008}.

In the past four decades, many tools have been developed for computing C*-envelopes,
and they can now be computed for many classes of examples.  
They have become a basic tool for the study of nonself-adjoint operator algebras.

Completely bounded maps took on a life of their own in the 1980s.
Wittstock \cite{Witt} proved the extension theorem for completely bounded maps.
Paulsen \cite{Paulsen_cpb} showed that this follows directly from 
Arveson's extension theorem using a $2\times2$ matrix trick. 

The notion of completely bounded and completely positive maps has had a profound 
influence on self-adjoint operator algebras.
While the goal of Bill's paper was dilation theory, the theory of completely
positive maps has had a life of its own in the self-adjoint theory. 
The injectivity of $\B(\H)$ led to deep work on injective von Neumann algebras by 
Connes \cite{Connes2}. 
This was followed by work on nuclear C*-algebras by Lance \cite{Lance},
and work by Choi, Effros and Lance determining the tight connection between 
nuclearity of a C*-algebra and the injectivity of
the von Neumann algebra it generates \cite{EffrosLance, ChoiEffros1, ChoiEffros3}.

When Brown, Douglas and Fillmore \cite{BDF1,BDF2} did their groundbreaking work 
on essentially normal operators and the Ext functor (a K-homology
theory for C*-algebras), Arveson pointed out how one obtains inverses in the Ext group using
completely positive maps. Later he wrote an important paper \cite{Arv1977} introducing
quasicentral approximate units for C*-algebras, and used this to provide a
unified and transparent approach to Voiculescu's celebrated generalized
Weyl-von Neumann theorem \cite{VoicWvN}, the Choi-Effros lifting theorem 
\cite{ChoiEffros2}, and the structure of the Ext groups for nuclear C*-algebras.

Completely positive and completely bounded maps also play a central role in
Kadison's similarity problem, which asks whether every bounded representation of
a C*-algebra is similar to a $*$-representation. This was confirmed for nuclear C*-algebras
by Bunce \cite{Bunce} and Christensen \cite{Christ}. 
Then Haagerup \cite{Haag_sim} established the result for cyclic representations 
of arbitrary C*-algebras. 
Along the way, he showed that a representation was similar to a $*$-representation 
if and only if it was completely bounded.
See Pisier's monograph \cite{Pisier_sim}.

The theory of CP and CB maps led to a considerable revolution on the 
nonself-adjoint side as well.
One important application concerns the Sz.Nagy-Halmos problem.
An operator $T$ is called \textit{polynomially bounded operator} 
if there is a constant $C$ so that
\[ \|p(T)\| \le C \|p\|_\infty \qforal p \in \bC[x] . \]
It is \textit{completely polynomially bounded} if this inequality holds 
for matrices of polynomials.
Halmos's refinement of Sz.Nagy's question asks whether 
every polynomially bounded operator is similar to a contraction.
Paulsen \cite{Paulsen_sim} showed that a representation of a
unital operator algebra is similar to a completely contractive one
if and only if it is completely bounded.
He used this to show that every completely polynomially bounded operator
is similar to a contraction \cite{Paulsen_cpb} . 
Then Pisier \cite{Pisier_pb} showed that there are polynomially
bounded operators which are not similar to a contraction,
solving the Sz.Nagy-Halmos problem.
See \cite{Paul2, Pisier_sim} for good coverage of this material.

The notions of matrix norms and completely bounded maps led to an abstract
categorical approach to abstract operator systems, operator algebras and operator spaces.
This is significant because it frees the theory from reliance on spatial representations.
Choi and Effros \cite{ChoiEffros4} characterized abstract operator systems.
Ruan \cite{Ruan} established a GNS style representation theorem for
operator spaces; and Blecher, Ruan and Sinclair \cite{BRS} did the
same for unital operator algebras. 
See \cite{EffrosRuan, Pisier_op, BleLeM} for a treatment of this material.

\section{Maximal subdiagonal algebras} \label{S:subdiag}

In his early work, Arveson used nonself-adjoint operator algebras to find
non-commutative analogues of results in analytic functions theory, to 
study dynamics on topological spaces, to study invariant subspaces, and even to 
study the corona problem of Carleson. 
We will survey these results in the next few sections.

K.~Hoffman \cite{Hoff_log} defined a \textit{logmodular} function algebra
to be a function algebra $\A$ on a compact Hausdorff space $X$
such that the set $\{\log |h| : h \in \A^{-1}\}$ is norm dense in $\rC_\bR(X)$.
He showed that this class shared many properties that the algebra $H^\infty$
of bounded analytic functions on the disk enjoys. 
In fact $H^\infty$ satisfies the stronger property that every invertible real
function $f$ in $L^\infty$ has the form $\log|h|$ for an invertible $H^\infty$ function.
Indeed one can take $h = \exp(f+i\tilde f)$ where $\tilde f$ is the harmonic conjugate of $f$.
We say that $H^\infty$ is strongly logmodular, or that it satisfies factorization.

Srinivasan and Wang \cite{SW} defined a weak-$*$ closed subalgebra $\A$ of $L^\infty(m)$
to be a weak-$*$ Dirichlet algebra if $\A+\A^*$ is weak-$*$ dense in $L^\infty(m)$ and
the functional $\phi(f)=\int f \,dm$ is multiplicative. Hoffman and Rossi \cite{HR} show
that these algebras are logmodular. A variety of other properties such as factorization,
variants of Szeg\"o's theorem and Beurling's invariant subspace theorem are
eventually shown to be equivalent to this property.
See \cite{BleLab_surv} for a survey of this and the non-commutative generalizations.

In the paper \cite{Arv1967b} that came out of Arveson's thesis \cite{Arv1964}, 
he developed a non-commutative setting for logmodularity. 
This setting is a finite von Neumann algebra $M$ with a faithful normal state $\tau$.
A weak-$*$ closed subalgebra $\A$ is tracial if there is a normal expectation $\Phi$
of $M$ onto $\A\cap\A^*$ which is a homomorphism on $\A$: $\Phi(ab) = \Phi(a) \Phi(b)$ 
for all $a,b\in\A$, and preserves the trace: $\tau\Phi = \tau$.
A (finite) \textit{maximal subdiagonal algebra} is a tracial subalgebra of $M$ 
such that $\A+\A^*$ is weak-$*$ dense in $M$.

It isn't hard to see that logmodularity of a function algebra is equivalent to saying
that $\{ |h|^2 : h \in \A^{-1}\}$ is dense in $\rC_\bR(X)$, and that factorization
is the stronger property that every invertible real continuous function on $X$
is exactly $|h|^2$ for some invertible $h\in\A$.
A subalgebra $\A$ of a C*-algebra $B$ is called logmodular if $\{a^*a : a\in\A^{-1}\}$
is dense in $B_{sa}$, and has factorization if this set coincides with the set
of all invertible self-adjoint elements.
By the polar decomposition, this latter property is equivalent to saying that every
invertible element of $B$ factors as $b=ua$ for some $a\in\A$ and unitary $u\in B$.

Arveson \cite{Arv1967b} established the analogue of the Hoffman-Rossi result:
\begin{thm}[Arveson 1967]
Every maximal subdiagonal algebra satisfies factorization.
\end{thm}

An interesting class of examples is obtained as follows.
Start with a chain $\L$ of projections (in the usual order) in a finite von Neumann algebra
$M$  such that there is a normal expectation onto the commutant $\L'$.
Then 
\[ \Alg\L = \{ a \in M : ae = eae \qforal e \in \L \} \]
is a maximal subdiagonal algebra.
Arveson also provides interesting examples associated to ergodic group actions
on measure spaces.

There are a variety of other important properties that $H^\infty$ satisfies,
and surprisingly many of them have natural non-commutative analogues.
To describe them, we need to introduce the Fuglede-Kadison determinant \cite{FK}
which they define on $II_1$ factors, 
and Arveson extends to all finite von Neumann algebras.
Define $\Delta(a) = e^{\tau(\log |a|)}$ when $a$ is invertible, and in general
\[ \Delta(a) = \inf_{\ep > 0} e^{\tau(\log (|a|+\ep 1))} \]
This map is multiplicative, positive homogeneous, monotone, self-adjoint;
and it is norm continuous on the invertible elements.
Arveson makes this much more tractable by showing that
\[ \Delta(a) = \inf \tau(|a|b) : b \in M^{-1},\ b > 0,\ \tau(b) \ge 1 \} .\]
{}From this it follows that 
\[ \Delta(a) \le \Delta(\Phi(a)) \qforal a \ge 0 .\] 

Now if $h\in H^\infty$, then Jensen's inequality says that
\[ \log |h(0)| \le \frac1{2\pi} \int_0^{2\pi} \log |h(e^{it})| \,dt .\]
Jensen's formula states that this is an equality when $h$ is invertible.
We also mention Szeg\"o's theorem, which states that if $f$ is a positive
function in $L^1$, then
\[
 \inf_{h\in H^\infty,\ h(0)=1} \frac1{2\pi} \int_0^{2\pi} |h(e^{it})|^2 f(t)\,dt = 
 \exp\Big( \frac1{2\pi} \int_0^{2\pi} \log f(t)\,dt \Big) .
\]

Arveson defines the analogues in the non-commutative setting: a tracial algebra $\A$
satisfies Jensen's inequality if 
\[ \Delta(\Phi(a)) \le \Delta(a) \qforal a \in \A \]
and satisfies Jensen's formula if
\[ \Delta(\Phi(a)) = \Delta(a) \qforal a \in \A^{-1}. \]
Szeg\"o's formula becomes:
\[ \Delta(b) = \inf\{ \tau(b a^*a) : a \in \A,\ \Delta\Phi(a) \ge 1 \}  .\]
Arveson \cite{Arv1967b} established that these three properties are equivalent, 
and that they are satisfied in many classes of examples.
Labuschagne \cite{Lab} proved that these properties hold for all
maximal subdiagonal algebras. Then Blecher and Labus\-chagne \cite{BleLab_ncH}
showed that for the class of tracial algebras, these properties are equivalent to
logmodularity, and to factorization, and to the weak-$*$ Dirichlet property.
See \cite{BleLab_surv} for a survey of this material.

\section{Dynamical systems} \label{S:dynamics}

Arveson was interested in using operator algebras to obtain complete
invariants for a discrete dynamical system in both the topological
and the measure theoretic settings. 

In the measure setting, consider an ergodic transformation $\tau$ of a 
measure space $(X,m)$ where $m$ is a probability measure. 
Two transformations $\sigma$ and $\tau$ are conjugate if there is a measure
preserving automorphism $\Phi$ of $X$ so that $\sigma = \Phi \tau \Phi^{-1}$.
An important problem is to determine when two ergodic transformations are conjugate.

This is easily transferred to a statement about operator algebras.
Take $M = L^\infty(m)$ considered as acting on $L^2(m)$ by multiplication.
The measure preserving transformation $\tau$ induces an automorphism of $M$ by 
\[
 \alpha(M_f) = M_{f\circ \tau}
 \quad\text{so that}\quad 
 \alpha(\upchi_E) = \upchi_{\tau^{-1}(E)} .
\]
Ergodicity means that the only fixed points of $\alpha$ are the scalars.
If $\sigma$ induces the automorphism $\beta$, then conjugacy of $\sigma$ and $\tau$
converts to conjugacy of $\alpha$ and $\beta$ in the automorphism group of $M$.

Every automorphism $\alpha$ of $M$ determines a unique unitary map $U$ as follows.
Since $L^\infty(m)$ is a dense subset of $L^2(m)$, we set
\[ U (g) = \alpha(M_g) 1 \qforal g \in L^\infty(m) .\]
Then one readily checks that this extends to a unitary operator such that 
\[ \alpha(M_g) = U M_g U^* \qforal g \in L^\infty(m) .\]
The conjugacy of two transformations reduces to unitary equivalence of the two
automorphisms.

Arveson's idea was to construct an operator algebra $\A(\alpha)$ as the norm
closure of
\[ \A_0(\alpha) = \Big\{ \sum_{i=1}^k M_{f_i} U^i : f_i\in L^\infty(m),\ k \ge 0 \Big\} .\]
His main result in \cite{Arv1967a} is that conjugacy of two automorphisms 
$\alpha$ and $\beta$ is equivalent to 
the unitary equivalence of $\A(\alpha)$ and $\A(\beta)$.

His results actually got a lot more mileage in the topological setting, 
which he and Josephson developed in \cite{ArvJos}. 
Begin with a locally compact Hausdorff space $X$ and a proper map 
$\sigma$ of $X$ into itself.  As we shall see, construction of a nonself-adjoint
operator algebra does not require the map $\sigma$ to be invertible.
Two dynamical systems $(X,\sigma)$ and $(Y,\tau)$ are (topologically) 
conjugate if there is a homeomorphism $\gamma$ of $Y$ onto $X$ 
so that $\sigma = \gamma\tau\gamma^{-1}$.
Again the issue is to decide when two systems are conjugate.

The map $\sigma$ induces an endomorphism $\alpha$ of $\rC(X)$ by composition.
When $\sigma$ is a homeomorphism, this map is an automorphism.
So we consider a dynamical system as the pair $(\rC(X), \alpha)$.
Arveson and Josephson construct a nonself-adjoint operator algebra in an 
analogous way to the measure theoretic setting.  Their construction depended
on the existence of a quasi-invariant measure and small fixed point set. 
Indeed they put rather stringent conditions on the maps:

Suppose that $(X,\sigma)$ and $(Y,\tau)$ are two topological dynamical systems
satisfying the following conditions
\begin{enumerate}
\item $\sigma$ and $\tau$ are homeomorphisms.
\item there is a regular Borel probability measure $\mu$ on $X$ which is quasi-inv\-ar\-iant for
$\sigma$, i.e. $\mu\circ\sigma$ and $\mu $ are mutually absolutely continuous;
\item $\mu(O)>0$ for every non-empty open set $O\subset X$;
\item the set $P = \{ x\in X : \sigma^n(x)=x \text{ for some }n\ge1 \}$ of periodic points
has measure zero;
\item there is a regular Borel probability measure $\nu$ on $Y$ with these
same properties, and in addition $\tau$ is ergodic with respect to $\nu$.
\end{enumerate}
Arveson and Josephson use these properties to build algebras 
$\A(\alpha)$ and $\A(\beta)$ as in the measure theoretic setting,
but using continuous functions rather than $L^\infty(\mu)$,
by constructing an explicit representation on $L^2(\mu)$.

What I will describe instead is a more general construction due to Peters \cite{Pet} 
called the \textit{semicrossed product}.
This has the major advantages that it does not depend on any measure,
and we may consider non-invertible systems $(X,\tau)$ where $\tau$ is any 
proper continuous map from $X$ into itself.

The key new notion is that of a \textit{covariant representation} $(\pi,V)$
where $\pi$ is a $*$-representation of $\rC(X)$ on a Hilbert space $\H$ and
$V\in\B(\H)$ is a contraction such that
\[  V\pi(f) = \pi(\alpha(f)) V \qforal f \in \rC(X) .\]
Peters defines the operator algebra $\rC(X) \times_\alpha \bZ_+$ to be the
universal operator algebra for covariant representations.
To form this algebra, let $\A_0(\alpha)$ be the set of all formal polynomials
of the form  $\sum_{i=0}^k f_i \tau^i$ for $f_i\in\rC(X)$ and $k\ge0$.
Define a norm by
\[
 \| \sum_{i=0}^k f_i \tau^i \| := 
 \sup \Big\{ \big\| \sum_{i=1}^k \pi(f_i)V^i  \big\| : (\pi,V) \text{ is covariant} \Big\} .
\]
The semicrossed product is the completion of $\A_0(\alpha)$ in this norm. 

Peters shows that a sufficient family of covariant representations is given
by the orbit representations: fix $x\in X$ and define a representation $\pi_{x}$
\[
 \pi_x(f) = \diag \big( f(x), f(\tau(x)), f(\tau^2(x)), \dots \big)
\]
and set $V$ to be the unilateral shift. 
The direct sum over a dense subset of $X$ yields a faithful completely isometric
representation of $\A(\alpha)$.  In particular, $\tau$ is represented by an isometry.

Peters shows that subject to the Arveson-Josephson hypotheses, 
their algebra $\A(\alpha)$ coincides with the semicrossed product
$\rC(X) \times_\alpha \bZ_+$.

I can now state the Arveson-Josephson result in this modern terminology.

\begin{thm}[Arveson-Josephson 1969]
Suppose that $(X,\sigma)$ and $(Y,\tau)$ are two topological dynamical systems,
and let $\alpha$ and $\beta$ be the automorphisms of $\rC(X)$ and $\rC(Y)$
induced by $\sigma$ and $\tau$ respectively.
Furthermore assume conditions (i)--(v) listed above.
Then the following are equivalent:
\begin{enumerate}
\renewcommand{\labelenumi}{(\arabic{enumi}) }
\item $(X,\sigma)$ and $(Y,\tau)$ are conjugate.
\item $\rC(X) \times_\alpha \bZ_+$ and $\rC(Y) \times_\beta \bZ_+$
are completely isometrically isomorphic.
\item $\rC(X) \times_\alpha \bZ_+$ and $\rC(Y) \times_\beta \bZ_+$
are algebraically isomorphic.
\end{enumerate}
\end{thm}

Peters improved on this result considerably by removing the 
conditions (i), (ii), (iii) and (v) at the expense of a slight strengthening of (iv) to 
\begin{enumerate}
\item[(a)] $X$ is compact.
\item[(b)] $\sigma$ has no periodic points. 
\end{enumerate}
He obtained the same equivalence.
Hadwin and Hoover \cite{HH} further improved on this.
Finally Davidson and Katsoulis \cite{DK_conj} removed all of these conditions,
and showed that (1), (2) and (3) are equivalent for arbitrary proper maps of locally compact
Hausdorff spaces with no conditions on periodic points at all.

\section{Invariant subspaces} \label{S:invariant}

The invariant subspace problem, which remains open today, was a popular
topic four decades ago. Arveson made some important contributions, which
as usual, stressed operator algebras rather than single operator theory.

If $\A$ is an algebra of operators on a Hilbert space $H$, let
\[ \Lat\A = \{ M \text{ is a closed subspace of }\H : \A M\subset M \} ,\]
and if $\L$ is a collection of subspaces, then
\[ \Alg(\L) = \{ T \in \B(\H): T M\subset M \foral M\in\L \} .\]
Halmos called an operator algebra \textit{reflexive} if $\A = \Alg\Lat\A$.
This algebra is always unital and closed in the weak operator topology (\wot-closed).
An operator algebra is \textit{transitive} if it has no proper closed invariant subspace.
The transitive algebra problem asks whether there is a proper unital
\wot-closed transitive operator algebra.

Arveson's first contribution to invariant subspace theory recognized the 
advantage of having a masa (maximal abelian self-adjoint algebra) in $\A$.
Note that he shows that the algebra is dense in the weak-$*$ topology, 
rather than weak operator topology.
This is a subtle, and occasionally non-trivial, strengthening.

\begin{thm}[Arveson's density theorem 1967]
 A transitive subalgebra of $\B(\H)$ which contains a masa is weak-$*$ dense in $\B(\H)$.
\end{thm}

This result received a fair bit of interest, and a generalization was found by
Radjavi and Rosenthal \cite{RR}: if $\A$ is an operator algebra
containing a masa, and $\Lat\A$ is a chain (nest), then 
the \wot-closure of $\A$ is reflexive.

Perhaps this result prompted Bill to revisit the issue. He decided to study
weak-$*$ closed operator algebras containing a masa, to determine whether all
such operator algebras are reflexive. We assume that the Hilbert space is separable.
The masa $\D$ is spatially isomorphic to $L^\infty(\mu)$ acting on $L^2(\mu)$
for some sigma-finite Borel measure $\mu$.
If $\A \supset \D$, then $\Lat\A \subset \Lat\D$.
If $M\in\Lat\D$, then the orthogonal projection $P_M$ commutes with $\D$;
and by maximality, $P_M$ belongs to $\D$.
So $\Lat\A$ can be naturally identified with a sublattice of the lattice of
projections in $\D$. In particular, these projections all commute.
In addition, the lattice is complete (as a lattice, and as a \sot-closed set of projections).
Such a lattice is called a \textit{commutative subspace lattice} (CSL).

In a 100 page tour de force in the Annals of Mathematics, 
Arveson \cite{Arv1974a} developed an extensive
theory for these operator algebras. The first major result is a spectral theorem for CSLs. 

\begin{thm}[Arveson 1974]
If $\L$ is a CSL on a separable Hilbert space, there is a compact metric space $X$, 
a regular Borel measure $\mu$ on $X$, and a Borel partial order $\preceq$ on $X$ 
so that $\L$ is unitarily equivalent to the lattice of projections
onto the (essentially) increasing Borel sets.
\end{thm}

Then he develops a class of pseudo-integral operators and identifies when one
is supported on the graph of the partial order; and hence belongs to $\Alg\L$.
Curiously, while the pseudo-integral operators seem to depend on the choice
of coordinates, the weak-$*$ closure of the pseudo-integral operators supported
on the graph of the partial order is independent of all choices. He calls this
algebra $\A_{\text{min}}(\L)$.  He shows that 
$\Lat\A_{\text{min}}(\L) = \L$ and hence $\Lat\Alg\L=\L$
for all CSLs. One of the main results is:

\begin{thm}[Arveson 1974]
If $\L$ is a CSL on a separable Hilbert space, 
then there is a weak-$*$ closed operator algebra $\A_{min}(\L)$ containing a 
masa with $\Lat(\A_{min}(\L)) = \L$ so that: 
a weak-$*$ closed operator algebra $\A$ containing a masa 
has $\Lat(\A)=\L$  if and only if 
$ \A_{min} \subseteq \A \subseteq \Alg(\L) .$
\end{thm}

Surprisingly, it is possible for $\A_{\text{min}}(\L)$ to be properly contained in $\Alg(\L)$,
even if it is \wot-closed.
An example is constructed based on the failure of spectral synthesis on $S^3$
in harmonic analysis. 
He calls a lattice \textit{synthetic} if $\A_{\text{min}}(\L) = \Alg(\L)$.
Many examples of synthetic lattices are produced, including lattices which are 
generated by finitely many commuting chains, and lattices associated to certain 
ordered groups. This is a definitive answer to the question.

The area of CSL algebras has received a lot of attention.  See \cite{Arv1984}
for a survey of some of this material.  Also see \cite{DavNest} for further material.

\section{Nest algebras and the corona theorem} \label{S:nest}

Kadison and Singer \cite{KS} initiated a theory of triangular operator algebras.
While their algebras were not necessarily \wot-closed or even norm-closed,
Ring\-rose \cite{Ring1,Ring2} developed a related class that are \wot-closed
and reflexive. He calls a complete chain of subspaces $\N$ a \textit{nest}.
A \textit{nest algebra} is $\T(\N) := \Alg\N$, the algebra of all operators which are
upper triangular with respect to the nest $\N$. In the first paper, he 
characterized the Jacobson radical. In the second, he showed that any algebra
isomorphism between two nest algebras is given by a similarity.

Arveson's paper  \cite{Arv1975} was very influential  for the subject of nest algebras.
In particular, he developed a distance formula from an arbitrary 
operator to $\T(\N)$.
Note that if $N\in\N$ and $T\in\T(\N)$, the invariance of $N$ may be expressed
algebraically by $TP_N = P_N TP_N$, or equivalently $P_N^\perp TP_N = 0$.
Therefore if $A\in\B(\H)$ and $T\in\T(\N)$, we see that
\[ 
 \|P_N^\perp A P_N\| = \|P_N^\perp (A-T) P_N\| 
 \le \inf_{T\in\T(\N)} \|T-A\| = \dist(A,\T(\N)) .
\]
This inequality remains valid if we take the supremum over all $N \in \N$.
Surprisingly, this leads to an exact formula:

\begin{thm}[Arveson's distance formula 1975]
Let $\N$ a nest on a Hilbert space $\H$; and let $A \in \B(\H)$.
Then 
\[ \dist(A, \T(\N)) = \sup_{N\in\N} \| P_N^\perp A P_N\| . \]
\end{thm}

\noindent
This and other results about nest algebras proved important in the development 
of the theory. 

A separable continuous nest is a nest which is order isomorphic to $I=[0,1]$. 
The Volterra nest $\N$ consisting of subspaces 
\[ N_t=L^2(0,t) \subset \H=L^2(0,1) \qfor 0\le t \le 1 \]
is an example which is multiplicity free, in that the projections
onto $N_t$ generate a masa. 
The infinite ampliation $\M$ given by 
$M_t = L^2((0,t)\times I) \subset L^2(I^2)$ generates
an abelian von Neumann algebra with non-abelian commutant.
So $\T(\N)$ and $\T(\M)$ are not unitarily equivalent.
Ringrose \cite{Ring2} asked whether they could be similar.

Three students of Arveson's eventually solved this problem completely.
Andersen \cite{Ander} showed that any two continuous nests are
approximately unitarily equivalent. So in particular, there is a sequence of unitaries
$U_n$ so that 
\[ [0,1] \ni t \to P_{M_t} - U_nP_{N_t} U_n^* \]
are continuous, compact operator valued functions converging uniformly to 0.
(Arveson \cite{Arv1983} generalized this to certain CSL algebras.)
Larson \cite{Larson} used this to solve Ringrose's problem by showing that
any two continuous nests are similar.
Shortly afterwards, Davidson \cite{Dav_sim} showed that if two nests are
order isomorphic via a map which preserves dimensions of differences of
nested subspaces, then there is a similarity of the two nests implementing
the isomorphism; and the similarity may be taken to be a small compact
perturbation of a unitary operator. 
See \cite{DavNest} for more material on nest algebras.

Arveson was also interested in finding an operator theoretic proof of the famous 
corona theorem of Carleson \cite{Carl}. This result is formulated as follows:
suppose that $f_1,\dots,f_n$ are functions in $H^\infty$ for which there is
an $\ep>0$ so that
\begin{align*} \tag{\dag}
 \inf_{z\in\bD}\  \sum_{i=1}^n |f_i(z)|^2 &\ge \ep .
\end{align*}
Then there are functions $g_i\in H^\infty$ so that $\sum_{i=1}^n f_ig_i = 1$.
A reformulation of this result is that the unit disk $\bD$ is dense in the maximal ideal
space of $H^\infty$.  
Carleson's proof was very difficult, and people continue to look for more accessible proofs.
An easier (but not easy) proof, based on unpublished ideas of Thomas Wolff, 
is now available.  See Garnett's book \cite{Garn}.

Arveson had an operator theoretic approach to the corona theorem.
$H^2$ is the Hardy space of square integrable analytic functions on the disk $\bD$
considered as a subspace of $L^2(\bT)$. 
Represent $H^\infty$ on $H^2$ as multiplication operators $T_h$ 
(these are analytic Toeplitz operators).  
In \cite{Arv1975}, he established this variation of Carleson's theorem: 

\begin{thm}[Arveson's Toeplitz corona theorem 1975]
If $f_1,\dots,f_n \in H^\infty$ and 
\begin{align*} \tag{\ddag}
 \sum_{i=1}^n T_{f_i}T_{f_i}^* \ge \ep I ,
\end{align*}
then there are $g_i \in H^\infty$ so that $\sum_{i=1}^n f_i g_i = 1$.
\end{thm}

This result is weaker than the corona theorem because $(\ddag)$ is a
stronger hypothesis than $(\dag)$.
To see this, we use the fact that $H^2$ is a reproducing kernel Hilbert space:  the functions 
\[ k_w(z) = \frac{(1-|w|^2)^{1/2}}{1-z\bar{w}} \]
are unit vectors for $w\in\bD$  with the property that $T_h^*k_w = \ol{h(w)} k_w$.  
Applying the vector state $\rho_w(A) = \ip{A k_w,k_w}$ to $(\ddag)$ for $w \in \bD$,
one recovers Carleson's condition $(\dag)$.

It turns out that Toeplitz corona theorems are useful in establishing full corona theorems.
Ball, Trent and Vinnikov \cite{BTV} established a Toeplitz corona theorem for all 
complete Nevanlinna-Pick kernels as a consequence of their commutant lifting theorem.
This was a useful step towards a bona fide corona theorem by Costea, Sawyer and Wick
\cite{CSW} for Drury-Arveson space, among other spaces.
For definitions of these concepts, see section~\ref{S:multivariable}.

\section{Automorphism groups} \label{S:auto}

In \cite{Arv1974b}, Arveson developed a spectral analysis of automorphism
groups acting on C*-algebras and von Neumann algebras.

Here $G$ is a locally compact abelian group with Haar measure $dt$.
A representation is a  bounded representation $t \to U_t$ in $\B(X)$,
where $X$ is a Banach space, and the representation is assumed to be continuous 
with respect to a certain weak topology. 
In particular, \sot-continuous unitary representations on Hilbert space form an 
important special case.
A standard technique extends this to a bounded representation $\pi_U$ 
of $L^1(G)$ such that
\[ \pi_U(f)\xi = \int_G f(t) U_t\xi \,dt .\]

We may regard the Fourier transform as the Gelfand transform of the commutative Banach
algebra $L^1(G)$ into $\rC(\hat G)$, where $\hat G$ is the dual group consisting of 
continuous characters of $G$.
When $f\in L^1(G)$, its zero set is $Z(f) = \{ \gamma \in \hat G : \hat f(\gamma) = 0 \}$.

We define the spectrum of the representation $U$ to be
\[ \Sp U = \bigcap\{ Z(f) : f \in L^1(G),\ \pi_U(f) = 0 \} .\]
Moreover, if $\xi\in \H$, then the spectrum of $\xi$ is
\[ \Sp_U(\xi) = \bigcap\{ Z(f) : f \in L^1(G),\ \pi_U(f) \xi = 0 \} .\]
Then if $E$ is a closed subset of $\hat G$, we define the \textit{spectral subspace} to be
\[ M^U(E) = \{ \xi : \Sp_U(\xi) \subset E \} .\]
When $U$ is a strongly continuous unitary representation on Hilbert space, 
Stone's theorem yields a projection valued measure $dP_\gamma$ on $\hat G$ so that 
$U_t = \int_{\hat G} \gamma(t) \,dP_\gamma$; and $M^U(E) = P(E)\H$.
But in more general contexts, there is often no version of Stone's theorem.
Nevertheless, the spectral subspaces form a substitute, and they determine the
representation uniquely:

\begin{thm} [Arveson 1974]
Let $X$ be a Banach space, and let $U$ and $V$ be two strongly continuous
representations of $G$ on $X$.  If $M^U(E)=M^V(E)$ for every compact subset $E$
of $\hat G$, then $U=V$.
\end{thm}

An important application of this theory is to a one parameter group of automorphisms 
acting on a von Neumann algebra $\R \subset \B(\H)$.  
That is, $t\to\alpha_t$ is a weak-$*$ continuous representation of $\bR$ into $\Aut(\R)$.
Arveson was interested in finding conditions that determine when this group is unitarily 
implemented, i.e.\ when $\alpha_t = \Ad U_t$, where $U$ is a strongly continuous unitary 
representation with $U_t\in\R$ and $\Sp U \subset [0,\infty)$.
Let $\R^\alpha(E) = \{ A \in \R : \Sp_\alpha(A) \subset E \}$.
The answer is given by:

\begin{thm} [Arveson 1974]
Let $(\alpha_t)$ be a weak-$*$ continuous representation of $\bR$ into $\Aut(\R)$.
Then there is a strongly continuous unitary representation $U$
with $U_t\in\R$ and $\Sp U \subset [0,\infty)$ such that $\alpha_t = \Ad U_t$ for $t\in\bR$ 
if and only if
\[\ \bigcap_{t \in\bR} \R^\alpha[t,\infty) = \{0\} .\]
Moreover $U_t = \int_{-\infty}^{+\infty} e^{itx} \,dP(x)$ where $P$ is the spectral measure
such that 
\[ P[t,\infty)\H = \bigcap_{s<t}  \R^\alpha[s,\infty)\H .\]
\end{thm}

Arveson used this to provide a new proof of the famous theorem of Kadison \cite{Kad} 
and Sakai \cite{Sakai} that every derivation of a von Neumann algebra is inner.
Alain Connes \cite{Connes1} used the Arveson spectrum to define the more refined
Connes spectrum of an automorphism group that allows one to define the type 
$III_\lambda$ factors.

One can show that 
\[ \R^\alpha[s,\infty) \R^\alpha[t,\infty)  \subset \R^\alpha[s+t,\infty) .\]
It follows that $H^\infty_\alpha(\R) = \R^\alpha[0,\infty)$ is a weak-$*$ closed subalgebra of $\R$.
This is sometimes a maximal subdiagonal algebra, but not always.
A very curious result of Loebl and Muhly \cite{LoeblMuhly} exhibits an example where
$H^\infty_\alpha(\R)$ is a proper nonself-adjoint subalgebra of $\R$ which is 
\wot-dense in $\R$.
A variant on the transitive algebra problem deals with \textit{reductive} algebras,
meaning that every invariant subspace is reducing (or that the orthogonal complement of
every invariant subspace is invariant).  
It asks whether every \wot-closed reductive operator is a von Neumann algebra.
Loebl and Muhly's construct is a reductive \textit{weak-$*$ closed} algebra which is 
not self-adjoint (but the \wot-closure is a von Neumann algebra).

\bigskip
I will now skip ahead to some more recent results. As I mentioned in the introduction,
in the intervening period, Bill was very interested in $E_0$-semigroups and CP-semigroups.
This work will be covered in the article \cite{Izumi} by Izumi in this volume.

\section{Multivariable operator theory} \label{S:multivariable}

As discussed in section \ref{S:dilation}, the Sz.Nagy dilation theorem led to a
development of function-theoretic techniques in single operator theory.
In particular, if $T$ is a contraction with no unitary summand, then it has
an $H^\infty$ functional calculus. This has many ramifications.  
See \cite{SF} for an updated look at this material.

The examples of Parrott and Varopolous, also mentioned in section \ref{S:dilation},
explain why there are problems with a straightforward generalization of this
dilation theory to several variables by studying $d$-tuples of commuting contractions.

It turns out that a different choice of norm condition leads to a better theory.
The more tractable norm condition on an $d$-tuple $T_1,\dots,T_d$ is a
row contractive condition: consider $T=\big[ T_1,\ \dots,\ T_d \big]$ as an
operator mapping the direct sum $\H^{(d)}$ of $d$ copies of $\H$ into $\H$
and require that
\[ \|T\| = \big\| \sum_{i=1}^d T_iT_i^* \big\|^{1/2} \le 1 .\]

In the non-commutative case, the appropriate model for the dilation theorem is 
best described as the left regular representation of the free semigroup $\Fd$,
consisting of all words in an alphabet $\{1,\dots,d\}$, acting on $\Fock$ by
\[ L_i \xi_w=\xi_{iw} \qfor 1 \le i \le d,\ w \in \Fd .\]
The space $\Fock$ can be realized as the full Fock space over $E=\bC^d$,
namely
\[ \Fock \simeq \sum_{k\ge0}^\oplus E^{\otimes k} = \bC \oplus \bC^d \oplus \bC^{d^2} \oplus \dots \]
where $E^{\otimes k} =\bC^{d^k}$ is identified with the span of words
of length $k$. In this view, the operators $L_i$ become the \textit{creation operators}
that tensor on the left by $e_i$, where $e_1,\dots,e_d$ is the standard basis for $\bC^d$.

The dilation theorem is that every row contraction $T$ dilates to a row isometry
$V = \big[ V_1,\ \dots,\ V_d \big]$ where $V_i ^* V_j= \delta_{ij} I$.
Furthermore, there is a decomposition $V_i \simeq L_i^{(\alpha)} \oplus W_i$
where $\alpha$ is a cardinal and $W_1,\dots,W_d$ are generators of a representation 
of the Cuntz algebra $\O_d$.
This was established for $d=2$ by Frazho \cite{Frazho}, for finite $d$ by Bunce \cite{Bun}
and for arbitrary $d$ including infinite values, and uniqueness, by Popescu \cite{Pop_diln}.
Popescu has written a long series of papers developing the Sz.Nagy-Foia\c s\ theory
in this context. The \wot-closed algebras generated by row isometries have also
been extensively studied by Popescu beginning with \cite{Pop_mult} and
Davidson and Pitts beginning with \cite{DP1}.

The dilation theory for commuting row contractions came earlier, but was exploited later.
Drury \cite{Drury} showed that every commuting row contraction with $\|T\|<1$
dilates to a direct sum of copies of a certain $d$-tuple of weighted shifts, 
$S= \big[ S_1,\ \dots,\ S_d \big]$, now called the \textit{$d$-shift}.
This was refined by M\"uller and Vasilescu \cite{MV} to show that if $\|T\|\le1$, 
then $T$ dilates to a row contraction $V$ with commuting entries so that
$V_i \simeq S_i^{(\alpha)} \oplus N_i$ where the $N_i$ are commuting normal operators
such that $\sum_{i=1}^d N_iN_i^* = I$, a \textit{spherical isometry}.
Perhaps the reason this did not go further at the time is that the exact nature of
the $d$-tuple $S$ was not realized.

Arveson began a program of multivariable operator theory in \cite{Arv1998}.
He reproved the Drury-M\"uller-Vascilescu dilation theorem. 
Moreover, he showed that the $d$-shift has a natural representation on symmetric
Fock space over $E=\bC^d$.  Let $E^k$ denote the symmetric tensor product
of $k$ copies of $E$, which can be considered as the subspace of $E^{\otimes k}$ 
which is fixed by the action of the permutation group $S_k$ which permutes the
terms of a tensor product of $k$ vectors.  Then
\[  H^2_d = \sum_{k\ge0} \!\strut^\oplus\, E^k \subset \Fock .\]
This space is coinvariant for the creation operators $L_i$, and $S_i = P_{H^2_d}L_i|_{H^2_d}$.

There is a natural basis 
\[ z^k = \frac{k!}{|k|!} \sum_{|w|=k,\ w(z) = z^k} \xi_w , \]  
where $k=(k_1,\dots,k_d) \in \bN_0^d$, $|k| = \sum_i k_i$, $k! := \prod_i k_i!$
and $z^k := z_1^{k_1}z_2^{k_2} \cdots z_d^{k_d}$.
This is an orthogonal but not orthonormal basis since $\|z^k\|^2 = \tfrac{k!}{|k|!} =: c_k$.
In this basis, $S_i$ is just multiplication by $z_i$.
The space $H^2_d$ may now be considered as the space of analytic functions
\[
 f(z) = \sum_{k\in \bN_0^d} a_k z^k 
 \quad\text{such that}\quad  
 \sum_{k\in \bN_0^d} |a_k|^2 c_k < \infty .
\]
These series converge on the unit ball $\bB^d$ of $\bC^d$, so are bona fide functions.
This is a reproducing kernel Hilbert space (RKHS), meaning that $w\in\bB^d$ separate points
and the value of $f(w)$ may be recovered by
\[ f(w) = \ip{f,k_w} \quad\text{where}\quad k_w(z) = \frac1{1-\ip{z,w}},\ w\in\bB^d .\]
The space $H^2_d$ is now often called Drury-Arveson space.

Every RKHS has an algebra of multipliers, those functions which multiply the space
into itself. 
The multiplier algebra $\M_d=\Mult(H^2_d)$ of $H^2_d$ consists of bounded 
analytic functions on the ball $\bB_d$. 
It is the \wot-closed unital algebra generated by $S_1$, \dots, $S_d$.
However the norm is greater than the sup norm, and is not comparable---so $\M_d$
is a proper subalgebra of $H^\infty(\bB_d)$.

The C*-algebra $\ca(S):=\ca(S_1,\dots,S_n)$ contains the compact operators $\fK$,
and $\ca(S)/\fK \simeq \rC(\partial \bB^d)$, the boundary sphere of $\bB^d$---which is
the spectrum of the generic spherical isometry. So one deduces that this is the
C*-envelope of the unital nonself-adjoint algebra $\A_d$ generated by $S_1$, \dots, $S_d$.

The classical Nevan\-linna-Pick theorem \cite{Pick} for the unit disk states that
given $z_1,\dots,z_n\in\bD$ and $w_1,\dots,w_n\in\bC$, there is a
function $h\in H^\infty$ with $\|h\|_\infty \le 1$ so that $h(z_i)=w_i$ if and only if
\[
  \Bigg[ \frac{1-w_i \ol{w_j}}{1-z_i \ol{z_j}} \Bigg]_{n\times n}  
\]
is positive definite.  A matrix version states that given $z_1,...,z_n$ in the disk, 
and $r \times r$ matrices $W_1$, \dots, $W_n$, there is a function $F$ in the the unit
ball of $\M_r(H^\infty)$ such that $F(z_i) = W_i$ if and only if the Pick matrix
\begin{align*}
 \Bigg[ \frac{I_r-W_iW_j^*}{1-z_i\ol{z_j}} \Bigg]_{n\times n}
\end{align*}
is positive semidefinite.  Sarason \cite{Sar} put this into an operator theoretic context.

An RKHS $\H$ of functions on a set $X$ is a natural place to study interpolation. 
One asks the same question: given $x_1,\dots,x_n \in X$ and  
$r \times r$ matrices $W_1$, \dots, $W_n$, is there a function $F$ in the the unit
ball of $\M_r(\Mult(\H))$ so that $F(z_i) = W_i$? It is easy to show that a necessary
condition is the positive semidefiniteness of the matrix
\begin{align*}
 \Big[ (I_r-W_iW_j^*)\ip{k_{x_j},k_{x_i}} \Big]_{n\times n} .
\end{align*}
A space for which this is also a sufficient condition is called a complete
Nevan\-linna-Pick kernel.

Complete Nevanlinna-Pick kernels were characterized by Quiggin \cite{Quig}
and McCullough \cite{McCull}. It follows from their characterization that $H^2_d$
is a complete Nevanlinna-Pick kernel. This was shown from a completely different
angle, by showing that the non-commutative \wot-closed algebra $\fL_d$ on Fock 
space generated by $L_1,\dots,L_d$ has a good distance formula for ideals
\cite{DP3, AriasPop},
and hence has a Nevanlinna-Pick theory based on Sarason's approach.
It follows that $\M_d$ is a complete quotient of $\fL_d$, and hence is a
complete Nevanlinna-Pick kernel.
Finally, Agler and McCarthy \cite{AgMcC} give a new proof of the Quiggin-McCullough 
theorem, and show that $\M_d$ for $d=\infty$ is the universal complete 
Nevanlinna-Pick kernel in the sense that every irreducible complete Nevanlinna-Pick RKHS
imbeds in a natural way into $H^2_d$.
See \cite{AgMcCbook} for an overview of this theory.

The confluence of these very different directions leading to the same space, $H^2_d$,
is remarkable: the good dilation theory for commuting row contractions,
the close relationship with the non-commutative model $\fL_d$,
and the fact that this is a complete Nevanlinna-Pick kernel.
All play a role in making $H^2_d$ a very fertile environment in which to study 
multivariable operator theory. 

Associated to any row contraction $T$ is the completely positive contractive map
\[ \phi(A) = \sum_{i=1}^d T_iAT_i^* .\]
Notice that $T$ is spherical if and only if $\phi(I)=I$.  
We call a row contraction $T$ \textit{pure} if there is no spherical part.
The decreasing sequence
\[ I \ge \phi(I) \ge \phi^2(I) \ge \dots \]
always has a limit in the strong operator topology.
This limit $\phi_\infty(I)$ is $0$ if and only if $T$ is pure.

In \cite{Arv2000}, Arveson begins a much more algebraic approach to multivariable
operator theory.  Douglas particularly has espoused the approach of considering
a Hilbert space $\H$ with an operator algebra $\A$ acting on it via a representation
$\rho:\A\to\B(\H)$ as a Hilbert module over $\A$ with the action $a\cdot\xi = \rho(a)\xi$.
This approach was developed for function algebras in the monograph \cite{DougPaul}
by Douglas and Paulsen. 

In this context, one begins with a commuting row contractive $d$-tuple $T$ acting on $\H$, 
and consider $\H$ as a module over $\bC[z] :=\bC[z_1,\dots,z_d]$ given by
$p\cdot\xi = p(T_1,\dots,T_d)\xi$. In order to use results from algebra effectively,
assume that the defect operator
\[ \Delta := (I - TT^*)^{1/2} =  \big(I -  \sum_{i=1}^d T_iT_i^* \big)^{1/2} \]
has \textit{finite rank}; and define $\rank(T) := \rank(\Delta)$.
We also say that $\H$ is a Hilbert module of the same rank.

Form a finitely generated module over $\bC[z]$ as follows
\[ M_\H := \{ p\cdot \xi : p \in \bC[z],\ \xi \in \Delta\H \} .\]
This space is \textit{not closed}; so it is a finitely generated $\bC[z]$ module algebraically.
{}From commutative algebra (Hilbert's syzygy theorem), one obtains a free resolution
\[ 0 \to F_n \to F_{n-1} \to \dots \to F_1 \to M_\H \to 0 .\]
Each $F_i$ is a finitely generated free $\bC[z]$-module of rank $\beta_i < \infty$.
One can define an invariant known as the Euler characteristic
\[ \upchi(\H) := \sum_{i=1}^n (-1)^{i+1}\beta_i \]
which is independent of the choice of resolution.

Motivated by an analogy with the Gauss-Bonnet theorem relating the curvature of a manifold to its topological invariants, Arveson also defines an analytic invariant that
he calls \textit{curvature}. Set 
\[ T(z) = \sum_{i=1}^d \ol{z}_i T_i = [ T_1\ \dots\ T_d ] [z_1\ \dots \ z_d]^* .\] 
This satisfies $\|T(z)\| \le \|z\|_2$.  Thus we may define a function on the ball by
\[
 F(z) = \Delta (I-T(z)^*)^{-1} (I-T(z))^{-1} \Delta|_{\Delta\H}
 \in \B(\Delta\H) \qfor z \in \bB_d .
\]
The boundary values, appropriately weighted, exist in the sense that
\[
 K_0(\H) = \lim_{r\to 1^-} (1-r^2) \Tr F(r\zeta) 
 \quad\text{exists a.e. for }\zeta\in \partial\bB_d .
\]
Define the \textit{curvature} of $\H$ to be
\[ K(\H) = \int_{\partial\bB_d} K_0(\zeta) \,d\sigma(\zeta) \]
where $\sigma$ is normalized Lebesgue measure on the sphere $\partial\bB_d$.

\begin{thm}[Arveson 2000]
Let $\H$ be a finite rank Hilbert module. Then
\[
 K(\H) = d! \lim_{n\to\infty} \frac{\Tr(I-\phi^{n+1}(I))}{n^d}
\]
and
\[
 \upchi(\H) = d! \lim_{n\to\infty} \frac{\rank(I-\phi^{n+1}(I))}{n^d} .
\]
\end{thm}

It is immediate that $0 \le K(\H) \le \upchi(\H)$.  
Examples show that this inequality may be proper.
However there is an important context in which they are equal.
Say that $\H$ is \textit{graded} if there is a strongly continuous unitary
representation $\Gamma$ of the circle $\bT$ on $\H$ so that 
\[
 \Gamma(t) T_i \Gamma(t)^* = t T_i 
 \qforal 1 \le i \le d,\ t \in \bT .
\]
This is called a gauge group on $\H$.
When $\H$ is graded, Fourier series allows the decomposition of $\H$ into
a direct sum of subspaces $\H_n = \{ \xi \in \H : \Gamma(t)\xi = t^n \xi \}$.
When $\H$ is finite rank, these are all finite dimensional.
Arveson calls the following result his analogue of the Gauss-Bonnet theorem:

\begin{thm}[Arveson 2000]
Let $\H$ be a graded finite rank Hilbert module. Then 
\[K(\H) =  \upchi(\H).\]
\end{thm}

In particular, the curvature is always an integer in the graded case.
Arveson hypothesized that it was always an integer. This conjecture
was verified by Greene, Richter and Sundberg \cite{GRS}.

An important tool for multivariable operator theory is the Taylor spectrum
and functional calculus \cite{Taylor1, Taylor2}.
Arveson \cite{Arv2002b} builds a Dirac operator based on this theory.
Let $Z=\bC^d$ with basis $e_1,\dots,e_d$;
and let $\Lambda Z$ be the exterior algebra over $Z$, namely
\[ \Lambda Z = \Lambda^0 Z \oplus \Lambda^1 Z \oplus \cdots \oplus \Lambda^d Z \]
where $\Lambda^k Z$ is spanned by vectors of the form $x_1\wedge \dots \wedge x_k$.
For $z\in Z$, there is a Clifford operator given by
\[ C(z) (x_1\wedge \dots \wedge x_k) = z \wedge x_1\wedge \dots \wedge x_k .\]
Now form $\tilde\H = \H \otimes\Lambda Z$, and for $\lambda \in \bC^d$ and $z \in Z$, define
\[
 B_{\lambda} = \sum_{i=1}^d (T_i-\lambda_i I) \otimes C(e_i),\quad 
 D_{\lambda} = B_{\lambda} + B_{\lambda}^* \qand
 R(z) = I \otimes (C(z)+C(z)^*) .
\]
Arveson observes that $\lambda$ is in the Taylor spectrum of $T$ if and only if
$D_{\lambda}$ is invertible. The pair $(D_0,R)$ is called the Dirac operator for $T$.

Now the graded space $\widetilde\H$ splits into its even and odd parts:
\[
 \widetilde\H_+ = \sum_{i \text{ even}} \oplus (\H \otimes \Lambda^i Z)
 \qand 
 \widetilde\H_- = \sum_{i \text{ odd}} \oplus (\H \otimes \Lambda^i Z).
\]
Note that $D_0$ maps $\widetilde\H_+$ into $\widetilde\H_-$ and vice versa.
Let $D_+ = D_0|_{\widetilde\H_+}$ considered as an operator in $\B(\widetilde\H_+, \widetilde\H_-)$.
The relationship to curvature is the following index theorem:

\begin{thm}[Arveson 2002]
If $T$ is a pure graded row contraction, then both $\ker D_+$ and $\ker D_+^*$
are finite dimensional, and 
\[ (-1)^d K(\H) = \ind D_+ = \dim\ker D_+ - \dim\ker D_+^* .\]
\end{thm}

Stability of the Fredholm index leads to the corollary that two $d$-contractions
which are unitarily equivalent modulo the compact operators have the same curvature.

Finally we discuss \cite{Arv2005} in which Arveson considers pure graded
finite rank $d$-contractions.
By the dilation theorem for commuting row contractions, 
such a row contraction is obtained from a graded
submodule $M$ of $H^2_d \otimes \bC^n$, where $n$ is finite. 
Graded means that $M$ is spanned by homogeneous polynomials, and thus
is generated by a finite set of homogeneous polynomials (Hilbert's basis theorem).

An easy calculation shows that the multipliers $S_i$ on $H^2_d$ satisfy
\[ [S_i, S_j^*] = S_iS_j^*-S_j^*S_i \in \fS_p \qforal 1 \le i,j \le d,\ p>d ,\]
where $\fS_p$ is the Schatten $p$-class of compact operators with 
$s$-numbers lying in $\ell^p$. 
Examples suggest that this condition on the commutators $[T_i,T_j^*]$
should persist for graded submodules $M$, i.e.\ $T_i = S_i^{(n)}|_M$.
Arveson establishes this for modules generated by monomials 
of the form $z_1^{k_1}\dots z_d^{k_d}\otimes \xi$. 
Moreover he makes the following conjecture:

\begin{conj}[Arveson 2005]
If $T=[T_1,\dots,T_d]$ is a pure graded finite rank $d$-contraction, then
$[T_i,T_j^*] \in \fS_p$ for all $1 \le i,j \le d$ and $p>d$.
\end{conj}

Douglas \cite{Douglas2006} further refines this conjecture first by enlarging the family of
Hilbert modules considered, and more significantly by considering the ideal
$I$ of all polynomials annihilating the $d$-tuple $T$. Let $Z$ be the zero set
of this ideal. Then Douglas conjectures that the commutators $[T_i,T_j^*]$
should lie in $\fS_p$ for all $p> \dim(Z \cap \bB_d)$.

There are no known counterexamples, so both of these conjectures remain open.
There has been a lot of recent interest.
The best result so far, due to Guo and Wang \cite{GuoWang}, establishes Arveson's
conjecture when $d=2$ or $3$, and when $M$ is singly generated.

\bigskip
\textit{Acknowledgements.} 
I thank Matthew Kennedy, Vern Paulsen, Gilles Pisier and David Pitts  for reading
the first draft of this paper and making many helpful suggestions that improved 
the final version.


\end{document}